\documentclass[12pt]{amsart}
\usepackage{amsmath}\usepackage{amsrefs}
\usepackage{amsfonts}\usepackage{verbatim}
\usepackage{amssymb}\usepackage{marvosym}
\usepackage{dictsym}
\usepackage{amstext}
\usepackage{amsbsy}
\usepackage{amsopn}\usepackage{MnSymbol}
\usepackage{amsthm}\usepackage{color}

\theoremstyle{definition}

\theoremstyle{remark}

\def\proclaim#1{\vskip0.5em\noindent{\bf #1}\it}
\def\endproclaim{\vskip0.5em\par\noindent\rm}
\def\proclaim#1{\vskip0.5em\noindent{\bf #1}\it}
\def\endproclaim{\vskip0.5em\par\noindent\rm}
\def\demo#1{\vskip0.5em\noindent{\bf #1\ }}

\def\cite#1{[#1]}
\def\text#1{\mbox{#1}}
\def\flushpar{\par\noindent}

\def\mod{\mbox{ mod }}
\newcommand{\mapright}[1]{%
    \smash{\mathop{%
        \hbox to 1cm{\rightarrowfill}
        }
    \limits^{#1}
    }
}
\newcommand{\mapleft}[1]{%
    \smash{\mathop{%
        \hbox to 1cm{\rightarrowfill}
        }
    \limits_{#1}
    }
}

\def\e{\epsilon}
\def\a{\alpha}

\def\G{\Gamma}

\def\d{\delta}
\def\D{\Delta}
\def\s{\sigma}

\def\th{\theta}
\def\l{\lambda}
\def\x{\times}

\def\f{\flushpar}

\def\v{\varphi}

\def\om{\omega}
\def\Om{\Omega}
\def\B{\mathcal B}

\def\({\biggl(}
\def\){\biggr)}

\def\<{\langle}
\def\>{\rangle}
\def\bul{\smallskip\f$\bullet\ \ \ $}
\def\Smi{\smallskip\f{ \Smiley\ \ \ }}
\def\lfl{\lfloor}\def\rfl{\rfloor}\def\sms{\smallskip\f}
\def\sms{\smallskip\f}\def\sbul{\f$\bullet\
\ \ $}\def\sms{\smallskip\f}\def\Smi{\smallskip\f{\Smiley\ \ \ }}
\def\lra{\longrightarrow}\def\Lra{\Longrightarrow}
\def\lfl{\lfloor}\def\rfl{\rfloor}\def\lcl{\lceil}\def\rcl{\rceil}\def\xbm{(X,\B,m)}\def\incss{\Bbb N^\Bbb N(\uparrow)}\def\xyr{\xrightarrow}

\begin{document}\title{  IP-rigidity and eigenvalue groups }
\author{Jon. Aaronson,\ Maryam Hosseini        $\&$  Mariusz Lema\'nczyk}
\address[Aaronson]{\ \ School of Math. Sciences, Tel Aviv University,
69978 Tel Aviv, Israel.}
\email{aaro@tau.ac.il}\address[Hosseini]{Faculty of Mathematical Science,
  University of Guilan, Rasht, Iran}
\email{ hoseini@guilan.ac.ir}\address[Lema\'nczyk]{Faculty
of Mathematics and Computer Science, Nicolaus Copernicus
University, ul. Chopina 12/18, 87-100 Toru\'n,
Poland}
\email{mlem@mat.uni.torun.pl}\date{May 24th 2012}

\begin{abstract}
We examine the class of increasing sequences of natural numbers which are IP-rigidity sequences for some weakly mixing probability preserving transformation. This property is closely related to the uncountability of the eigenvalue group of  a corresponding non-singular transformation. We give examples, including a  super-lacunary sequence which is not IP-rigid.
\end{abstract}
\subjclass[2010]{37A05, (37A30, 37A40)}\keywords{Probability preserving transformation, spectral type, rigidity, IP-convergence, Dirichlet measure,  eigenvalue group.}

\thanks{Aaronson and Hosseini would like to thank   Nicolaus Copernicus University, and Aaronson the University of Surrey, for hospitality provided when this paper was done. Aaronson's research was partially
supported by Israel Science Foundation grant No. 1114/08. }
\maketitle\markboth{Jon. Aaronson,\ Maryam Hosseini        $\&$  Mariusz Lemanczyk \copyright  2012}{IP-rigidity and eigenvalue groups}
\section*{\S0 Introduction}\
Let $(X,\B,m)$ be a standard, continuous, probability space (that is, $(X,\B)$ is a Polish space equipped with its Borel sets and a non-atomic $m\in\mathcal P(X)$ (the collection of probability measures on $(X,\B)$).
 \

 We'll denote by ${\tt MPT}={\tt MPT}(X,\B,m)$ the collection of invertible, probability preserving transformations of $(X,\B,m)$. This is a Polish space when equipped with the {\it coarse topology} with basic neighborhoods of form
\begin{align*}U(T_0,f_1,&\dots,f_N;\e):=\\ &\{T\in\text{\tt MPT}:\ \|f_j\circ T^s-f_j\circ T_0^s\|_{L^2(m)}<\e\ \forall\ 1\le j\le N,\ s=\pm 1\}\end{align*}
where $T_0\in\text{\tt MPT}\xbm$ and $f_1,\dots,f_N\in L^2(m)$.
\par Equipped with the coarse topology,  {\tt MPT}$\xbm$ is a topological group under composition. It is embedded into the Polish,  topological  group $\mathcal U(L^2(m))$ of unitary operators on $L^2(m)$ equipped with the strong operator topology by the {\it Koopman representation} $U_Tf:=f\circ T\ \ (T\in\text{\tt MPT}\xbm,\ f\in L^2(m))$.
Accounts of the spectral theory of  unitary   operators can be found in \cite{KT} and \cite{N2}.
\subsection*{Recurrence and rigidity}
\

A  sequence $q\in\Bbb N^\Bbb N(\uparrow):=\{q\in\Bbb N^\Bbb N:\ q_n<q_{n+1}\ \forall\ n\ge 1\}$  is called a sequence of {\em recurrence} for  $T\in{\tt MPT}(X,\B,m)$ if $\limsup_{n\to\infty}\mu(A\cap T^{-q_n}A)>0$ for each $A\in\B$ of positive measure.
\

Rigidity is a stronger version of recurrence.
\

An  sequence $q\in\Bbb N^\Bbb N(\uparrow)$ is called a {\em rigidity sequence} for $T\in{\tt MPT}\xbm$ if $\mu(T^{q_n}A\triangle A)\to 0$ for each $A\in\B$; equivalently
\begin{align*}&\tag{\dsmilitary}
T^{q_n}\xrightarrow[n\to\infty]{\tt MPT} \text{\tt Id}.\end{align*}

 Using spectral theory one sees that (\dsmilitary) is equivalent to the {\tt restricted spectral type} $\s_T$ of $T$ (i.e. $U_T|_{L^2(m)_0}$) having the {\it Dirichlet property} along $q$, that is
\begin{align*}\tag{\ddag}
\chi_{q_n}\xrightarrow[n\to\infty]{L^2(\Bbb T,\sigma_T)}\,1.\end{align*} where $\chi_k(t):=e^{2\pi i kt}$.

\

Rigidity sequences for non-trivial  transformations  must be sparse. In particular,  unless $T\in{\tt MPT}$ is purely periodic any rigidity sequence for $T$ has at most finite intersection with each of its translates whence has zero  Banach density.

Additional properties of  rigidity  sequences are studied in \cite{BJLR} $\&$ \cite{EG} including the rigidity properties of {\tt lacunary} and {\tt super-lacunary} sequences, a sequence $q\in\incss$ being called
 {\it lacunary} if $\tfrac{q_{n+1}}{q_n}\ge\l>1\ \forall\ n\ge 1$ and  {\it super-lacunary} if $\tfrac{q_{n+1}}{q_n}\xyr[n\to\infty]{}\infty$.
\

\subsection*{Rigid factors, mild mixing and IP sets}\ \ \ \ Let $T\in{\tt MPT}\xbm$ and let $q\in\incss$. It is well known that the collection of sets
\begin{align*}\tag{\Wheelchair}\mathcal R(q):=\{A\in\B:\ m(A\D T^{q_n}A)\xyr[n\to\infty]{}0\}\end{align*}
is a $T$-invariant, $\s$-algebra. It corresponds to the maximal factor of $T$ which is rigid along $q$.
The transformation $T\in{\tt MPT}\xbm$ is called {\it mildly mixing} if it has no non-trivial,  rigid factor along any $q\in\incss$ (as in \cite{FW}).
\

Since the spectral type $\s_S$ of a factor $S$ of $T$ is absolutely continuous with respect to $\s_T$,
it is evident that $T$ has some non-trivial rigid factor if and only if $\exists$ a  {\it Dirichlet measure} $\mu\ll\s_T$,
(that is, one satisfying (\ddag) along some $q\in\incss$).
\

An {\it IP-set} is a  collection of ``finite sum sets" of form
$$\text{\tt FS}\,(q):=\{q(F):\ F\in\mathcal F\}$$ where $q\in\incss$ and for
$$F\in\mathcal F:=\{\text{\tt finite, nonempty subsets of}\ \Bbb N\},\ q(F):=\sum_{j\in F}q_j.$$
This notion appears in combinatorics, ultrafilter theory, topological dynamics (see \cite{Fu} and \cite{HS}) and also in ergodic theory.
\

As shown in  \cite{Fu}, $T\in{\tt MPT}\xbm$ is mildly mixing if and only if $\exists\ K\subset\Bbb N$ which intersects with every finite sum set so that
$$m(A\cap T^{-n}B)\xyr[n\to\infty,\ n\in K]{}m(A)m(B)\ \forall\ A,\ B\in\B;$$
equivalently (see  \cite{HMP1}),
\

$T$ is not mildly mixing if and only if
 $\exists\ q\in\incss$ so that
$$\underset{n\in\text{\tt\tiny FS(q)}}\inf\,|\widehat{\mu}(n)|\ >0.$$

The considerations involved   give rise to the notion of
 \subsection*{ IP convergence}
 \

 Let $q\in\incss$. We'll say that a sequence $a:\Bbb N\to Z$ (a metric space) {\it converges IP to $L\in Z$ along {\tt FS}$(q)$}
$$\text{\rm written}\ \ a(n)\xrightarrow[n\to\infty]{\text{\tt\tiny FS}\,(q)} \ L,\ \text{\rm in}\ Z\ \ \ \text{\rm  if}$$
$$a(q(F))\xrightarrow[F\in\mathcal F,\ \min\,F\to\infty]{Z} \ L.$$
This paper is about
\subsection*{IP-rigidity}
\

 We'll say that $b\in\Bbb N^\Bbb N(\uparrow)$ is

an {\it IP-rigidity sequence for  $T$} and that $T$ is {\it IP-rigid along $b$} if $$T^{n}\xrightarrow[n\to\infty]{\text{\tt\tiny FS}\,(b)}\ \text{\tt Id\ {\rm in}\   MPT}.$$
\

Let
\begin{align*}&\text{\tt IPRWM}:=\{b\in\Bbb N^\Bbb N(\uparrow):\ \exists\ T\in\text{\tt MPT}\xbm,   \text{\rm\small weakly mixing $\&$ IP-rigid along}\ b\}.\end{align*}

Any rigid transformation is IP-rigid on some subsequence (see \cite{Fu}). On the other hand if a transformation is
IP-rigid on $q\in\incss$, then it is rigid along much thicker subsequences (see \S5).
\

Similarly to (\Wheelchair), for $q\in\incss$, the collection
$$\mathcal R_{\tiny\tt IP}(q):=\{A\in\B:\ m(A\D T^{n}A)\xyr[n\to\infty]{{\tt FS}\,(q)}\ 0\}$$
is a $T$-invariant, $\s$-algebra. It corresponds to the maximal factor of $T$ which is IP-rigid along $q$.
As above, $T$ has a non-trivial factor, IP-rigid along $q$ if and only if
$\lim_{N\to\infty}\inf_{F\in\mathcal F,\ \min\,F\ge N}|\widehat{\s_T}(q(F))|>0.$
\

The existence of IP-Dirichlet measures along $b\in\Bbb N^\Bbb N(\uparrow)$ is related to the groups
\begin{align*}&G_p(b):=\{t\in\Bbb T:\ \sum_{n=1}^\infty\|b_nt\|^p<\infty\}\  (0<p<\infty)\ \&\\ & G_\infty(b):=\{t\in\Bbb T:\ \|b_nt\|\xyr[n\to\infty]{}0\}\end{align*}
where for $x\in\Bbb R,\ \|x\|:=\min_{n\in\Bbb Z}\,|x-n|$.
\

These groups are discussed in [AN] and [HMP2].
\

\subsection*{\large Results}
\

\proclaim{Proposition 1}\ \ \ \ Suppose that
$b\in \Bbb N^\Bbb N(\uparrow)$, then
\begin{align*}\ \ \ \ |G_1(b)|>\mathbf{\aleph_0}\ \ \Lra\ \ b\in \text{\tt IPRWM}.\end{align*} \endproclaim
Proposition 1 (which is folklore) can be proved using Propositions  1.1 and 1.2 (below).
\

\proclaim{Theorem 2} If $b\in\incss$, then
\begin{align*} b\in \text{\tt IPRWM}\ \Lra\ \ |G_2(b)|>\mathbf{\aleph_0}.\end{align*}
\endproclaim

Theorem 2 also provides an answer to a question in \cite{BJLR}:
 {\Large\Smi} {\em if $b\in\text{\tt IPRWM}$, then some irrational rotation is rigid along $b$}
  because if $|G_2(b)|>\aleph_0$ then $\exists\ \a\in G_2(b)\setminus\Bbb Q$.
   It follows that rotation of $\Bbb T$ by $\a$ is rigid along $b$. \Checkedbox

The converse of theorem 2 holds for   {\tt arithmetic} sequences  and {\tt Erdos-Taylor} sequences for different reasons.
\

A sequence $q\in\incss$ is called {\it arithmetic} if it is either
\bul {\it multiplicative} in the sense that $q_n|q_{n+1}\ \forall\ n\ge 1$; or it is the
\bul {\it  principal denominator sequence} of some  $\a\in\Bbb T\setminus\Bbb Q$, being defined by
$q_0=1,\ q_1=a_1,\ q_{n+1}:=a_{n+1}q_n+q_{n-1}$ where $\a=[0;a_1,a_2,\dots]$ is the continued fraction expansion of $\a$.
\proclaim{Proposition 3} Let $b\in\incss$, be arithmetic. The following statements are equivalent.
\begin{align*}&\text{\rm (a)}\  \ \varlimsup_{n\to\infty}\frac{b_{n+1}}{b_n}=\infty.\ \ \text{\rm (b)}\   |G_1(b)|>\mathbf{\aleph_0}.\\ &
\ \text{\rm (c)}\ b\in \text{\tt IPRWM}.\ \ \text{\rm (d)}\ \ |G_2(b)|>\mathbf{\aleph_0}.\end{align*}
\endproclaim
\

The {\it Erdos-Taylor  sequence associated to $(a_1,a_2,\dots)\in\Bbb N^\Bbb N$} is
 $b=(b_1,b_2,\dots)\in\Bbb N^\Bbb N(\uparrow)$ defined by
 $$b_1:=1.\ \ b_{n+1}:=a_{n}b_n+1.$$
 Erdos-Taylor  sequences were introduced in \cite{ET} and are considered to be ``extremely non-arithmetic''.
\proclaim{Proposition 4} If $b\in\incss$, is an Erdos-Taylor sequence, then
\begin{align*}\text{\rm (i)}\   \sum_{n\ge 1}\(\frac{b_n}{b_{n+1}}\)^2<\infty\ \iff\ \text{\rm (ii)}\   b\in \text{\tt IPRWM}\ \iff\ \text{\rm (iii)}\  \ |G_2(b)|>\mathbf{\aleph_0}.\end{align*}
\endproclaim
 We'll see that there are super-lacunary Erdos-Taylor sequences $b\ \&\ q\in\incss$ satisfying
  $|G_2(b)|>\mathbf{\aleph_0}\ \&\ G_1(b)=\{0\}$  and  $G_2(q)=\{0\}$.

\subsection*{Eigenvalue Groups and theorem 2}
\

 Groups of form $G_2$ appear as eigenvalue groups (see \cite{AN}). Eigenvalue groups and rigidity are related as follows:
\

An ergodic probability preserving transformation $S$ is not mildly mixing (i.e. has a rigid factor) if and only if there is a conservative, ergodic non-singular transformation $T$ so that $S\x T$ is not ergodic (see \cite{FW}). By the ergodic multiplier theorem of M. Keane (see e.g. \S2.7 of \cite{A}), this situation is characterized by $\s_S(e(T))>0$ where $\s_S$ is the restricted spectral type of $S$ and $e(T)$ is the group of eigenvalues of $T$.
\

We prove Theorem 2 in \S4 by considering a  dyadic cocycle (see below) associated to $b\in\Bbb N^\Bbb N(\uparrow)$ over the dyadic adding machine.
The  eigenvalue group of the  Mackey range (as in p. 76-77 in \cite{Z}) of this cocycle   is $G_2(b)$.
In case $b$ is a {\it growth sequence} as in \cite{A2}, that is $b(n)>\sum_{1\le k<n}b(k)$, then the Mackey range preserves a $\s$-finite measure and is isomorphic to the appropriate dyadic tower over the dyadic adding machine (defined in \cite{A2}).
\subsection*{Organization of the paper}
\

In \S1 we establish the basic results on Dirichlet sets and measures and begin to consider membership of {\tt IPRWM}.
\

In \S2 we  consider the class of arithmetic sequences,  and prove Proposition 3.
\

In \S3 we prove proposition 4  for Erdos-Taylor sequences and give our   main  examples.
\

The proofs  of propositions  3 $\&$ 4  both use Theorem 2 which is established  in \S4. In \S5 we make some quantitative remarks on the growth of rigid sequences for transformations IP-rigid along some
(particular) $b\in\incss$.
\

\section*{\S1 Dirichlet sets and measures}
\subsection*{Dirichlet sets}

A {\it Dirichlet set} is a subset  $\G\subset\Bbb T$ of form
$$\G(b)=\{t\in\Bbb T:\ \chi_{b_n}(t)\xrightarrow[n\to\infty]{}\ 1\}$$
where $b\in\Bbb N^\Bbb N(\uparrow)\ \&\ \ \chi_n(t):=e^{2\pi int}$.\par An {\it IP-Dirichlet set} is a subset  $\G\subset\Bbb T$ of form
$$\G(\text{\tt FS}\,(b))=\{t\in\Bbb T:\ \chi_{n}(t)\xrightarrow[n\to\infty]{\text{\tt\tiny FS}\,(b)} \ \ 1\}$$ where $b\in\Bbb N^\Bbb N(\uparrow)$. Here, we have
\proclaim{Proposition 1.1}\ \ \ For $b\in\Bbb N^\Bbb N(\uparrow)$,
$$\G(\text{\tt FS}\,(b))=G_1(b).$$\endproclaim
\demo{Proof sketch of $\subseteq$}
\

 It suffices to show that for $t\in\Bbb R$,
$$\|nt\|\xrightarrow[n\to\infty]{\text{\tt\tiny FS}\,(b)} \ \ 0\ \Rightarrow\ \ \sum_{n\ge 1}\|b_nt\|<\infty.$$

For $x\in\Bbb R$, let $\lfl x\rcl$ be the nearest integer to $x$ (if there are two, take the lesser one), and let
$$\<x\>:=x-\lfl x\rcl,\ \text{then}\ \ \ |\<x\>|=\|x\|\le \frac12.$$
Fix $t\in\Bbb R$ so that $\|nt\|\xrightarrow[n\to\infty]{\text{\tt\tiny FS}\,(b)} \ \ 0$, let $K>0$ be so that
\begin{align*}\|b(F)t\|<\frac1{16}\ \ \ \forall\ F\in\mathcal F,\ \min\,F\ge K.\end{align*}
If $F,\ G\subset [K,\infty)\cap\Bbb N$ are disjoint finite sets, then
$$\|b(F)t\|,\ \|b(G)t\|,\ \|b(F\cupdot G)t\|<\frac1{16}.$$
Since $\<b(F\cupdot G)t\>-\<b(F)t\>-\<b(G)t\>\in\Bbb Z$, this forces
\f$\<b(F\cupdot G)t\>=\<b(F)t\>+\<b(G)t\>.$
\

It follows  that
$$\sum_{n\ge K,\ \<b_nt\>\ge 0}\|b_nt\|\le\frac1{16},\ \sum_{n\ge K,\ \ \<b_nt\><0}\|b_nt\|\le\frac1{16}\ \&\ \sum_{n\in\Bbb N,\ n\ge K}\|b_nt\|<\frac18.\ \ \ \text{\Checkedbox}.$$

\subsection*{Dirichlet measures}
\
 A probability measure $\mu\in\mathcal P(\Bbb T)$ is
called
\bul a  {\it  Dirichlet} measure if
$$\|\chi_{b_n}-1\|_{L^2(\mu)}\xrightarrow[n\to\infty]{} \ \ 0$$ for some   $b\in\Bbb N^\Bbb N(\uparrow)$ in which case $\mu$ is called {\it Dirichlet along $b$} and
\bul an {\it IP Dirichlet} measure if
$$\|\chi_n-1\|_{L^2(\mu)}\xrightarrow[n\to\infty]{\text{\tt\tiny FS}\,(b)} \ \ 0$$ for some    $b\in\Bbb N^\Bbb N(\uparrow)$ in which case $\mu$ is called {\it IP-Dirichlet along $b$}.

\

Evidently:
\sms $\chi_{n_k}\xrightarrow[k\to\infty]{L^2(\mu)} 1$ if and only if $\widehat\mu(n_k)\xrightarrow[k\to\infty]{} \mu(\Bbb T)$,
\sms if $\mu$ is IP-Dirichlet along $b$, then so is any $\nu\ll\mu$,
\sms if $\mu$ is IP-Dirichlet, then $\exists\ b\in\incss$ so that $\sum_{n\ge 1}\|\chi_{b_n}-1\|_{L^2(\mu)}<\infty$, whence $\mu(G_1(b))=1$ and $\mu$ is IP-rigid along $b$.
\

By Proposition 1.1, a totally atomic measure $\mu\in\mathcal P(\Bbb T)$ is IP-Dirichlet along $b$ if and only if $\mu(G_1(b))=1$.
Examples in \S4 (below) show that this is false for
 continuous measures $\mu\in\mathcal P(\Bbb T)$.
\

\

\proclaim{Proposition 1.2}
\begin{align*}\text{\tt IPRWM}\ =\ \{b\in\Bbb N^\Bbb N(\uparrow):\ \exists\ \mu\in\mathcal P(\Bbb T)\ \ \text{\rm\small continuous $\&$ IP-Dirichlet along}\ b\}.\end{align*}\endproclaim\demo{Proof of $\subseteq$}
Suppose that $b\in\text{\tt IPRWM}$ and that $(X,\B,m,T)$ is a weakly mixing, probability preserving transformation so that
$T^{n}\xrightarrow[n\to\infty]{\text{\tt\tiny FS}\,(b)}\ \text{\tt Id}$. Fix $f\in L^2(m),\ \int_Xfdm=0,\
\int_X|f|^2dm=1$. The spectral measure of $f$: $\mu\in\mathcal P(\Bbb T)$ is continuous and IP-Dirichlet along $b$.
\demo{Proof of $\supseteq$}\ \ Suppose that $\mu_0\in\mathcal P(\Bbb T)$ is continuous and IP-Dirichlet along $b$.
Let $\mu$ be the symmetrization of $\mu_0$ (also continuous and IP-Dirichlet along $b$) and let $(X,\B,m,T)$ be the shift of the Gaussian process
with spectral measure $\mu$. The  spectral type of $(X,\B,m,T)$ is $\s_T=\sum_{n\ge 0}\mu^{n*}$ where $\mu^{n*}$ denotes the $n$-fold convolution
of $\mu$ with itself (see e.g. \cite{CFS}). Each $\mu^{K*}$ is continuous (whence $T$ is weakly mixing) and IP-Dirichlet along $b$
(since $\widehat{\mu^{K*}}(n)=\widehat{\mu}(n)^{K}\xrightarrow[n\to\infty]{\text{\tt\tiny FS}\,(b)}\ 1$). Every $\nu\ll\s_T$ is also IP-Dirichlet along $b$
and $T$ is IP-rigid along $b$. Thus $b\in\text{\tt IPRWM}$.\  \Checkedbox

\

It follows from proposition 1.2 that if $b\in\incss$ is an IP rigidity sequence for some $T\in\text{\tt MPT}$ not of discrete spectrum (i.e. $\s_T$ is not totally atomic) then $b\in\text{\tt IPRWM}$.

\

We complete this section with a ``mixed" multiplicative-finite sum condition for membership in {\tt IPRWM}.
\proclaim{Proposition 1.3}\ \ Suppose that $b\in\incss$ and that $\exists\ S\subset\Bbb N$ infinite, so that
$\sum_{n\in S}\frac{b_n}{b_{n+1}}<\infty$, and $b_n|b_{n+1}$ for $n\notin S$, then $b\in {\tt IPRWM}$.\endproclaim
If $\Bbb N\setminus S$ is finite, then $|G_1(b)|>\aleph_0$ by Theorem 5 in \cite{ET}, whence
$b\in {\tt IPRWM}$  by proposition 1.

 \demo{Proof} \ \ Assume (without loss of generality)  that $b_1=1$.    We  construct a weakly mixing
 $T\in {\tt MPT}\xbm$ which is IP-rigid along $b$ by cutting and stacking as in Ch. 7 of \cite{N1}.
\

To this end, we construct a nested sequence of Rokhlin towers $(\tau_n)_{n\ge 1}$ of intervals where $\tau_n$ has height $b_n$.

   Let $\tau_1$ be $[0,1]$. To construct $\tau_{n+1}$  from $\tau_n$:
    \bul if $n\notin S\ \&\ b_{n+1}=a_n b_n,\ a_n\in\Bbb N,\ a_n\ge 2$ then  we   cut $\tau_n$ into $a_n$ columns  and stack.
    \bul If $n\in S\ \&\ b_{n+1}= a_nb_n+r_n$, $ a_n,\ r_n\in\Bbb N,\ 1\le r_n< b_n$, we cut $\tau_n$ into $ a_n$ columns and put one spacer interval above the $\lfl\frac{ a_n}2\rfl$'th column from the left and  $r_n-1$ spacer intervals above the last column in the right and then we stack.
    \

    The tower $\tau_n$ is a stack of $b_n$ intervals of length $\prod_{j=1}^{n-1}\frac1{a_j}$ called {\it levels of $\tau_n$}.
    \

     It follows from \S7.31 in \cite{N1} that the transformation $T$ constructed preserves a finite measure $m$.
 A standard  argument as in the proof of Proposition 3.10 of \cite{BJLR} shows that $T$ is weakly mixing.

\

 Next,  we show that if $A$ is a union of levels in some $\tau_K$, then $$\sum_{n=1}^\infty m(A\D T^{b_n}A)<\infty.$$
 \

 To see this, we note first that $A$ is also a union of levels in every $\tau_n\ \ (n\ge K)$. Fix $n\ge K$ and write
 $$S\cap [K,\infty)=\{s_1<s_2<s_3<\dots\}.$$
 Since
 $\frac1{a_{s_\ell}+1}\le\ \frac{b_{s_\ell}}{b_{s_\ell+1}}\le \frac1{a_{s_\ell}}$ our assumptions imply
  $\sum_{\ell=1}^\infty \frac1{a_{s_\ell}}<\infty$.
 \

  To estimate $m(A\D T^{-b_n}A)$ for $n\in (s_{\ell-1},s_{\ell}]$ we consider the appearance of  the tower $\tau_n$ as ``$\tau_n$-stalks" inside $\tau_{s_\ell+1}$.

 A {\it $\tau_n$-stalk} in $\tau_{s_\ell+1}$ is a union  $\frak s=\bigcupdot_{k=0}^{b_n-1}T^kB$ of levels of $\tau_{s_\ell+1}$ where $B$ is contained in the base of $\tau_n$.

   Let $\frak s\in\B$ be a $\tau_n$-stalk  in $\tau_{s_\ell+1}$, then
$$m(\frak s\cap A)=\prod_{j=n}^{s_\ell}\frac1{a_j}\cdot m(A).$$
    \

  By construction, $\tau_{s_\ell+1}$ consists entirely of  $\tau_n$-stalks and spacer stalks  added in $\tau_{s_\ell+1}$. Thus, all  points of $A$ except   those   contained in the two  $\tau_n$-stalks preceding the   spacer stalks  added in $\tau_{s_\ell+1}$, return  to $A$ at  time  $b_n$, so
 $$ m(A\D T^{-b_n}A)=2m(A\setminus T^{-b_n}A)=4\prod_{j=n}^{s_\ell}\frac1{a_j}\cdot m(A)\le \frac4{2^{s_\ell-n}}\cdot\frac1{a_{s_\ell}}\cdot m(A).$$
 Thus, writing $s_0:=K-1$, we have
 \begin{align*}\sum_{n=K}^\infty  m(A\D T^{-b_n}A)&=
 \sum_{\ell=1}^\infty\sum_{n\in (s_{\ell-1},s_\ell]}m(T^{b_n}(A)\Delta A)\\ &\le
 \sum_{\ell=1}^\infty\sum_{n\in (s_{\ell-1},s_\ell]}\frac4{2^{s_\ell-n}}\cdot\frac1{a_{s_\ell}}\cdot m(A)\\ &\le 4m(A)
 \sum_{\ell=1}^\infty \frac1{a_{s_\ell}}\\ &<\infty.\end{align*}
 It follows from this that for $A$ a union of levels in some tower $\tau_n$ and $F=\{n_1<n_2<\dots<n_k\}\in\mathcal F$,
 \begin{align*}m(A\D T^{b(F)}A)&=m(A\D T^{\sum_{j=1}^kb_{n_j}}A)\\ &\le
 m(A\D T^{b_1}A)+m(T^{b_1}A\D T^{\sum_{j=1}^kb_{n_j}}A)\\ &=
 m(A\D T^{b_1}A)+m(A\D T^{\sum_{j=2}^kb_{n_j}}A)\\ &\le\\ &\vdots\\ &\le\sum_{j=1}^k m(A\D T^{b_j}A)\\ &\le\sum_{n=\min\,F}^\infty m(T^{b_n}(A)\Delta A)\\ &
 \xrightarrow[F\in\mathcal F,\ \min\,F\to\infty]{} 0.\end{align*}
 The collection of measurable sets $\mathcal C$ satisfying this last convergence is a $\s$-algebra and
 $$\mathcal C\supset\s\(\bigcup_{n\ge 1}\{\text{unions of levels in }\ \tau_n\}\)=\mathcal B.$$ Thus  $T$ is IP-rigid along $b$.\ \ \Checkedbox
\subsection*{Remark 1.4}\ \ \ The proof of Proposition 1.3 establishes the following proposition:
\sms {\em \ Suppose that $b\in\incss,\ p>0$ and that $\exists\ S\subset\Bbb N$ infinite, so that
$\sum_{n\in S}(\frac{b_n}{b_{n+1}})^p<\infty$, and $b_n|b_{n+1}$ for $n\notin S$, then  $\exists\ T\in{\tt MPT}\xbm$ weakly mixing and a dense collection $\mathcal A\subset\B$ so that \begin{align*}\tag{\dscommercial$(p)$}\ \ \ \ \ \ \ \ \ \sum_{n\ge 1}m(A\D T^{b_n}A)^p<\infty\ \ \forall\ A\in\mathcal A.\end{align*}}
As  above, {\dscommercial$(1)$ $\Lra$ IP-rigidity along $b$ whence $\s_T(G_2(b)^c)=0$ by Theorem 2.
Using the spectral theorem, one sees that {\dscommercial$(\tfrac12)$ $\Lra$ $\s_T(G_1(b)^c)=0$.

\section*{\S2 Arithmetic sequences}
In this section, we  prove proposition 3. The implications  (b) $\Lra$ (c)\ $\Lra$ (d)  $\Lra$ (a) follow from proposition 1,  theorem 2 and theorem 16 in \cite{E} (respectively). None of these uses arithmeticity.

We turn to the remaining implication (a) $\Lra$ (b).

\proclaim{Lemma 2.1}
\

If $b\in\Bbb N^\Bbb N(\uparrow)$ is multiplicative and $\sup_{n\ge 1}\frac{b_{n+1}}{b_n}=\infty$, then
$|G_1(b)|>\mathbf{\aleph_0}$.\endproclaim C.f. theorem 3 in \cite{ET}.\demo{Proof}
Suppose that $b_{n+1}=a_{n+1}b_n$, where $a_n\ge 2\ \forall n\geq 1$.
\

Since $\sup_{n\ge 1}\,a_n=\infty$, $\exists$ a subsequence $(n_k)_{k\geq 1}$ such that  $a_{n_{k+1}}/a_{n_k}\geq 3\ \forall\ k\ge 1$. Define $t:\Om:=\{0,1\}^\Bbb N\to [0,1]$ by
$$t(\om):=\sum_{k=1}^\infty\frac{\om_k}{b_{n_k}}.$$ Since
$$\sum_{k=L+1}^\infty\frac{\om_k}{b_{n_k}}\le \frac1{b_{n_L}}\sum_{j=1}^\infty\frac1{a_{n_L+1}a_{n_L+2}\cdots a_{n_j}}\le\frac1{b_{n_L}}\sum_{j=1}^\infty\frac1{3^j}< \frac1{b_{n_L}},$$
we have that  $t:\Om:=\{0,1\}^\Bbb N\to [0,1]$ is strictly increasing (with respect to lexicographic order on $\Om$), whence injective and $|t(\Om)|>\mathbf{\aleph_0}$.
\

It suffices to show that $t(\Om)\subset G_1(b)$.
\

To see this, fix $\om\in\Om$. For $N\ge 1$, we have that
\begin{align*}
b_Nt(\om)=\sum_{k=1}^\infty\frac{b_N\om_k}{b_{n_k}}=
\sum_{k\ge 1,\ n_k\ge N}\frac{\om_k}{a_{N+1}a_{N+2}\cdots a_{n_k}}\ \ \mod\ 1,
\end{align*}
whence
\begin{align*}
\|b_Nt(\om)\|=
\sum_{k\ge 1,\ n_k\ge N}\frac{\om_k}{a_{N+1}a_{N+2}\cdots a_{n_k}}\le \sum_{k\ge 1,\ n_k\ge N}\frac1{a_{N+1}a_{N+2}\cdots a_{n_k}}=:\D_N.
\end{align*}
For $n_{K-1}<N\le n_K$,
\begin{align*}\D_N =
\sum_{k\ge K}\frac1{a_{N+1}a_{N+2}\cdots a_{n_k}}\le
\sum_{k\ge K}\frac1{a_{N+1}a_{N+2}\cdots a_{n_K}\cdot a_{n_k}}\le
\frac1{2^{n_K-N}a_{n_K}}.
\end{align*}
Thus

\begin{align*}\sum_{N=n_1}^\infty\|b_{N}t(\om)\|&\leq\sum_{N=n_1}^\infty\D_N=
\sum_{k=2}^\infty\sum_{N=n_{k-1}+1}^{n_k}\sum_{\nu=k}^\infty{1\over{a_{N+1}a_{N+2}\cdots a_{n_\nu}}}\\ &\le
\sum_{k=2}^\infty\sum_{N=n_{k-1}+1}^{n_k}\frac1{2^{n_k-N}a_{n_k}}\le
\sum_{k=2}^\infty\frac2{a_{n_k}}\le 4\end{align*}
and $t(\om)\in G_1(b)$.\ \ \ \Checkedbox

\proclaim{Lemma 2.2}
\

Let  $q=q(\a)\in\Bbb N^\Bbb N(\uparrow)$  be the principal denominator sequence of $\a\in (0,1)\setminus\Bbb Q$. \

If  $\sup_{n\ge 1}\frac{q_{n+1}}{q_n}=\infty$, then
$|G_1(q)|>\mathbf{\aleph_0}$.\endproclaim
\demo{Proof}
\ As shown in \cite{IN}, \par for any $t\in[0,1]$ there is a unique sequence $(\omega_n)_{n\geq 1}\in\prod_{n\geq 1}\{0,1,\cdots,a_n\}$ such that \sbul $\omega_k\leq a_k,\ \omega_k=a_k\ \Rightarrow\ \omega_{k+1}=0$ and \sbul $t=\sum_{n=1}^\infty \omega_n \<q_n\alpha\>\ \mod 1$.
\

Since $\sup_{n\ge 1}{a_n}=\sup_{n\ge 1}\frac{q_{n+1}}{q_n}=\infty$, we can choose a sub-sequence $a_{n_k}$  such that $\sum_{k\geq 1}{1/a_{n_k}}<\infty$ and define
$$t:\Om=\{0,1\}^\Bbb N\to \Bbb T\ \ \text{by}\ \ t(\om):=\sum_{k=1}^\infty \omega_{k} \<q_{n_k-1}\alpha\>\ \mod 1.$$
By the above, $t:\Om\to \Bbb T$ is injective.  It follows that $|t(\Om)|>\mathbf{\aleph_0}$ and  it suffices to show that
$t(\Om)\subset G_1(q)$.
\

We claim that $\sup_{\om\in\Omega}\sum_{n=1}^\infty \|q_nt(\om)\|<\infty$.

 \

 Fix $K\geq 1$ and consider $n_K-1\leq N\leq n_{K+1}-1$.
Then
\begin{align*}
\|q_Nt(\om)\|&\leq\sum_{k=1}^\infty \|\omega_{n_k-1}q_Nq_{n_k-1}\alpha\|\\
&\leq  \sum_{k=1}^\infty\|q_Nq_{n_k-1}\alpha\|\\
&\leq \sum_{k=1}^K{q_{n_k-1}\over q_{N+1}}+\sum_{k=K+1}^\infty{q_N\over{q_{n_k}}}.\\
\end{align*}
Using the fact that  $\frac{q_{j+n}}{q_j}\geq \sqrt{2}^{n-2}\ \forall\ j,\ n\ge 1$  we have for some absolute constant $C$,
$${q_{n_k-1}\over q_{N+1}}={q_{n_k-1}\over q_{n_K-1}}\cdot{q_{n_K-1}\over q_{n_K}}\cdot{q_{n_K}\over q_{N+1}}
\leq {C\over {\sqrt{2}^{n_K-n_k}}}\cdot{1\over a_{n_K}}\cdot {C\over{\sqrt{2}^{N-n_K}}}\ \ \ \text{for $k\leq K$}$$
and
$${q_N\over{q_{n_k}}}={q_N\over{q_{n_{K+1}-1}}}\cdot{q_{n_{K+1}-1}\over q_{n_{K+1}}}\cdot{q_{n_{K+1}}\over q_{n_k}}
\leq {C\over {\sqrt{2}^{n_{K+1}-N}}}\cdot{1\over{a_{n_{K+1}}}}\cdot{C\over \sqrt{2}^{n_k-n_{K+1}}}\ \ \ \text{for $k>K$.}$$
 Therefore,
 \begin{align*}
\sum_{N=n_K-1}^{n_{K+1}-2}\|q_Nt(\om)\|&\le \sum_{N=n_K-1}^{n_{K+1}-2}\left(\sum_{k=1}^K{q_{n_k-1}\over q_{N+1}}+\sum_{k=K+1}^\infty{q_N\over{q_{n_k}}}\right)\\
&\leq C^2\sum_{N=n_K-1}^{n_{K+1}-2}\left(\sum_{k=1}^K {1 \over {a_{n_{K}}\sqrt{2}^{N-n_k}}}+\sum_{k=K+1}^\infty {1\over a_{n_{K+1}}\sqrt{2}^{n_k-N}}\right)\\
&\leq  C^3\sum_{N=n_K-1}^{n_{K+1}-2}(\frac1{a_{n_{K}}\sqrt{2}^{N-n_K}}+\frac1{a_{n_{K+1}}\sqrt{2}^{n_{K+1}-N}})\\
&\leq C^4({1\over a_{n_K}}+{1\over a_{n_{K+1}}})
\end{align*}
 and

\begin{align*}
\sum_{N=1}^\infty \|q_Nt(\om)\|\le C^4\sum_{K=1}^\infty({1\over a_{n_K}}+{1\over a_{n_{K+1}}})<\infty.\ \ \ \text{\Checkedbox}
\end{align*}
Hence $q\in{\tt IPRWM}. \ \ \ \text{\Checkedbox}$
The proof of Proposition 3 is now complete.\ \ \Checkedbox

\section*{\S3 Super-lacunary sequences}
\

Suppose that $b=(b_1,b_2,\dots)\in\Bbb N^\Bbb N(\uparrow)$ is super-lacunary, i.e. $\frac{b_{n+1}}{b_n}\xyr[n\to\infty]{}\infty$.
 \

 As in  theorem 17 in \cite{E},
we fix $N\ge 1$ with $\frac{b_{n+1}}{b_n}>10\ \forall\ n\ge N$ and let
$$E:=\bigcap_{n\ge N}E_n$$ where
$$E_n:=\{t\in [0,1]:\ \|b_nt\|\le\frac{4b_n}{b_{n+1}}\}.$$
Now
$$E_n\supseteq\bigcup_{k=1}^{b_n-1}I_{k,n}$$
where
$$I_{k,n}:=[\frac{k}{b_n}-\frac4{b_{n+1}},\frac{k}{b_n}+\frac4{b_{n+1}}].$$
For each $n\ge N$, the intervals $\{I_{k,n}:\ 1\le k<b_n\}$ are disjoint, and
each interval $I_{k,n}$ contains at least five disjoint intervals of form $I_{k',n+1}$. It follows that $E$ contains a  Cantor set and $|E|>\aleph_0$. Thus,  (c.f  of theorem 5 in \cite{ET})
\proclaim{Proposition 3.1}\ \ \ Suppose that $b\in\incss\ \&\ p>0$, then
$$\sum_{n\ge 1}\(\frac{b_n}{b_{n+1}}\)^p<\infty\ \ \ \Lra\ \ |G_p(b)|>\aleph_0.$$\endproclaim
\demo{Proof}\ \ Let $N\ge 1\ \&\ E$ be as above, then for $t\in E$,
$$\sum_{n\ge 1}\|b_nt\|^p\le N+4^p\sum_{n\ge N}(\tfrac{b_n}{b_{n+1}})^p<\infty.\ \ \ \ \text{\Checkedbox}$$

\

\proclaim{Proposition 3.2}\ \ \ Suppose that $b\in\incss$ and that
$$\sum_{n\ge 1}\(\frac{b_n}{b_{n+1}}\)^2<\infty,\ \ \text{then}\ \ \ b\in{\tt IPRWM}.$$\endproclaim
\demo{Proof} \ \ By Proposition 1.2, it suffices to construct a continuous probability in $\Bbb T$ which is IP-Dirichlet along $b$.
\

To this end, let $N\ge 1,\ E$ and  $\{I_{k,n}:\ 1\le k<b_n\}$ be as above.  As above, for each $n\ge N$, the intervals $\{I_{k,n}:\ 1\le k<b_n\}$ are disjoint, and
each interval $I_{k,n}$ contains at least five disjoint intervals of form $I_{k',n+1}$.
\

Thus we may choose
$$I_n(\om)=I_{k_n(\om),n}\ \ \ (n\ge 1,\ \om\in\{0,1\}^n)$$ so that
$$I_{n+1}(\om,\e)\subset I_n(\om)\ \ \forall\ n\ge 1,\ \om\in\{0,1\}^n\ \&\ \e=0,1$$
and

$X_{n+1}(\om,0)<X_n(\om)<X_{n+1}(\om,1)$ where $X_n(\om):=\frac{k_n(\om)}{b_n}$.

\

Next, for $\om\in\Om:=\{0,1\}^\Bbb N$,
$$X_n(\om_1,\dots,\om_n)\xyr[n\to\infty]{}X(\om)\ \text{where}\ \ \bigcap_{n\ge 1}I_n(\om_1,\dots,\om_n)=\{X(\om)\}$$ and
\begin{align*}b_nX(\om)&=b_nX_n(\om)+b_n(X_{n+1}(\om)-X_n(\om))+b_n(X(\om)-X_{n+1}(\om))\\ &=k_n(\om)+\xi_n(\om)+\th_n(\om)
\end{align*}

where
$$\xi_n(\om):=b_n(X_{n+1}(\om)-X_n(\om))\ \&\ $$
$$\th_n(\om):=b_n(X(\om)-X_{n+1}(\om)).$$ Note that  $|\th_n(\om)|\le\mathcal E_n:=\frac{4b_n}{b_{n+2}}$ and that by assumption, $\sum_{n\ge 1}\mathcal E_n<\infty$.

For $n\ge 1,\ \om\in\{0,1\}^n,\ \exists\ !\ p_{n,\om}:\{0,1\}\to (0,1)$ so that
$$p_{n,\om}(0)+p_{n,\om}(1)=1\ \&\ \ X_{n+1}(\om,0)p_{n,\om}(0)+X_{n+1}(\om,1)p_{n,\om}(1)=X_n(\om).$$
Define $P:\{\text{\tt cylinders}\}\to (0,1)$ by
$$P([a_1,a_2,\dots,a_n]):=\frac12\prod_{k=1}^{n-1}p_{k,(a_1,a_2,\dots,a_k)}(a_k).$$
It follows that $P$ is additive and by standard extension theory $\exists\ \Bbb P\in\mathcal P(\Om)$ extending $P$.

\

Denoting expectation with respect to $\Bbb P$ by $\Bbb E$ and writing $\om=(\om_1,\om_2,\dots,\om_n)\in\{0,1\}^n$, we have
$$\Bbb E(X_{n+1}\|\om_1,\om_2,\dots,\om_n)=X_{n+1}(\om,0)p_{n,\om}(0)+
X_{n+1}(\om,1)p_{n,\om}(1)=X_n(\om),$$ whence
$$\Bbb E(\xi_n\|\om_1,\om_2,\dots,\om_n)=b_n(\Bbb E(X_{n+1}\|\om_1,\om_2,\dots,\om_n)-X_n(\om))=0$$ and $\Bbb E(\xi_n)=0$.
\

For $n,\ k\ge 1$,
\begin{align*}\Bbb E(\xi_n\xi_{n+k})&=b_nb_{n+k}\Bbb E(\Bbb E(\xi_n\xi_{n+k}\|\om_1,\om_2,\dots,\om_{n+k})) \\ &=
b_nb_{n+k}\Bbb E((\xi_n(\om)\Bbb E\xi_{n+k}\|\om_1,\om_2,\dots,\om_{n+k}))=0\ \ \ \ \&\
\end{align*}
$$\Bbb E(\xi_n^2)=:\D_n\ \le\ \frac{16b_{n}^2}{b_{n+1}^2}.$$
By assumption $\sum_{n\ge 1}\D_n<\infty$, so
$\sum_{n\in K}\xi_n$ converges in $L^2(\Bbb P)$ for every $K\subset\Bbb N$ and
$$\Bbb E\((\sum_{n\in K}\xi_n)^2\)=\sum_{n\in K}E(\xi_n^2)=\sum_{n\in K}\D_n.$$

\

The measure $\mu:=\Bbb P\circ X^{-1}\in\mathcal P(\Bbb T)$ is continuous. We claim that it is IP-Dirichlet along $b$.

To check this, let $F\subset\Bbb N\cap [K,\infty)$ be finite and write $\Xi_F:=\sum_{n\in F}\xi_n$, then
\begin{align*}\|\chi_{b(F)}-1\|_{L^2(\mu)}& \le\|\<\tsum_{N\in F}b_Nt\>\|_{L^2(\mu)}\\ &=
\|\<\sum_{N\in F}(\xi_N+\th_N)\>\|_{L^2(\Bbb P)}\\ &
\le\|\<\Xi_F\>\|_{L^2(\Bbb P)}+\sum_{N\in F}\mathcal E_N.\end{align*}
Next,
\begin{align*}\|\<\Xi_F\>\|^2_{L^2(\Bbb P)}&=\Bbb E(\<\Xi_F\>^2)\\ &=
\Bbb E(1_{[|\Xi_F|\le\frac12]}\<\Xi_F\>^2)+\Bbb E(1_{[|\Xi_F|>\frac12]}\<\Xi_F\>^2)\\ &\le
\Bbb E(\Xi_F^2)+\frac14\Bbb P[|\Xi_F|>\frac12])\ \ \ \because\ \<\Xi_F\>^2\le\frac14,\\ &\le
2\Bbb E(\Xi_F^2)\ \ \ \text{\rm by Tchebychev's inequality}\\ &=2\sum_{N\in F}\Bbb E(\xi_N^2)\le 2\sum_{N=K}^\infty\D_N.\end{align*}
Thus,
\begin{align*}\|\chi_{n(F)}-1\|_{L^2(\mu)}& \le\|\<\Xi_F\>\|_{L^2(\Bbb P)}+\sum_{N=K}^\infty\mathcal E_N\\ &\le
\sqrt{2\sum_{N=K}^\infty\D_N}+\sum_{N=K}^\infty\mathcal E_N\\ &\xyr[K\to\infty]{}0\end{align*}
 proving that $\mu$ is IP-Dirichlet along $b$.\ \ \ \Checkedbox
\subsection*{Remark}\ \ \ The converses to propositions 3.1 $\&$ 3.2 are false.
It is easy to construct $b\in\incss$ multiplicative, super-lacunary so that $\sum_{n\ge 1}(\tfrac{b_n}{b_{n+1}})^p=\infty \ \ \ \forall\ p>0$.
By proposition 3, $|G_1(b)|>\aleph_0\ \&\ \ T\in\text{\tt IPRWM}$.
\subsection*{\large Erdos-Taylor  sequences $\&$  proposition 4}
\

We begin with a strong converse to Proposition 3.1  for Erdos-Taylor  sequences:
\proclaim{Proposition 3.3}\ \ \ Suppose that $b\in\incss$ is an Erdos-Taylor  sequence and let $p>0$, then
$$\sum_{n\ge 1}\(\frac{b_n}{b_{n+1}}\)^p=\infty\ \ \ \Lra\ \ G_p(b)=\{0\}.$$\endproclaim
This was stated in \cite{ET} for $p=1$ and $b$ the Erdos-Taylor  sequence  associated to $(2,3,\dots)$. See also Th\'eor\`eme 2 in \cite{P}.
\demo{Proof}\ \
  Let $b\in\incss$ be the Erdos-Taylor  sequence associated to $(a_1,a_2,\dots)\in\Bbb N^\Bbb N$ and let $t\in\Bbb R\setminus\Bbb Z$, then
$$\|b_nt\|<\frac{\|t\|}{2a_n}\ \Lra\ \|b_{n+1}t\|>\frac{\|t\|}{2}.$$
If $p>0$ and $\sum_{n\ge 1}(\frac{b_n}{b_{n+1}})^p=\infty$, then for  $t\in\Bbb R\setminus\Bbb Z$,
\bul either $\|b_nt\|\ge\tfrac{\|t\|}{2a_n}$ eventually and  $\sum_{n\ge 1}\|b_nt\|^p=\infty$, or
\bul $\|b_{n+1}t\|>\tfrac{\|t\|}{2}$ infinitely often and  $\sum_{n\ge 1}\|b_nt\|^p=\infty$.
\

Either way, $t\notin G_p(b)$.\ \ \Checkedbox

\demo{Proof of Proposition 4} \ \ The implications (i) $\Lra$ (ii)  $\Lra$ (iii)  $\Lra$ (i) follow from proposition 3.2, theorem 2 and proposition 3.3 (respectively).\ \ \ \Checkedbox
\

\subsection*{Examples}
\

\bul  If $b\in\incss$ is the Erdos-Taylor sequence associated to $(2,3,\dots)$, then
$$\text{\rm (a)}\ \ \ \ b\in\text{\tt IPRWM} \ \ \ \ \&\ \ \ \text{\rm (b)}\ \ \ \ G_1(b)=\{0\}.$$

\bul If $b\in\incss$ is the Erdos-Taylor sequence associated to $(a_2,a_3,\dots)$ where $a_n\to\infty$ and
$\sum_{n\ge 1}\frac1{a_n^2}=\infty$ (e.g. $a_n:=\lfl\sqrt n\rfl$), then $b$ is super-lacunary, $ G_2(b)=\{0\}$ and $b$ is not a sequence of IP-rigidity for any    probability preserving transformation other than the identity.

\section*{\S4 Proof of Theorem 2}
\

We prove Theorem 2 using  dyadic cocycles over the  dyadic odometer.
\

Let $\Om:=\{0,1\}^\Bbb N,$ and let $P\in\mathcal P(\Om)$ be symmetric product measure: $P=\prod(\frac12,\tfrac12)$, and let $\tau:\Om\to\Om$ be the dyadic odometer
defined by
$$\tau(1,\dots,1,0,\om_{\ell+1},\dots)=(0,\dots,0,1,\om_{\ell+1},\dots)$$
where $\ell=\ell(\om):=\min\,\{n\ge 1:\ \om_n=0\}$.
\

The {\it dyadic cocycle $\varphi:\Om\to\Bbb Z$ associated to $b\in\incss$} is  defined by
$$\varphi(\om):=b_{\ell(\om)}-\sum_{k=1}^{\ell(\om)-1}b_k$$ and its {\it skew product}

$\tau_\v:\Om\x\Bbb Z\to \Om\x\Bbb Z$ is  defined by $\tau_\v(x,n)=(\tau(x),n+\v(x))$.

Define $Q:\Om\x\Bbb Z\to \Om\x\Bbb Z$ by $Q(x,n):=(x,n+1)$ and fix $p\in\mathcal P(\Om\x\Bbb Z),\ p\sim P\x\#$ where $\#$ is counting measure on $\Bbb Z$.
\

There is (see \cite{Z}, pp. 76--77)
 an ergodic, non-singular transformation $(\frak X,\B(\frak X),\frak q,\frak T)$ of a standard probability space  and a map $\pi:\Om\x\Bbb Z\to \frak X$ so that
 $$p\circ\pi^{-1}=\frak q,\ \ \pi^{-1}\B(\frak X)=\{A\in\B(\Om\x\Bbb Z):\ \tau_\v A=A\}\ \ \&\ \ \pi\circ Q=\frak T\circ\pi.$$
 The ergodic, non-singular transformation $(\frak X,\B(\frak X),\frak q,\frak T)$ is called the {\it Mackey range of $(\tau,\v)$.}
 \

 In case $b$ is a growth sequence, equivalently $\v:\Om\to\Bbb N$, there is a $\s$-finite, invariant $T$-invariant measure $\frak m\sim \frak q$ with respect to which the  Mackey range $(\frak X,\B(\frak X),\frak m,\frak T)$ is isomorphic to the tower over $(\Om,\B(\Om),P,\tau)$ with height function $\v$ (aka the {\it dyadic tower with growth sequence $b$} in \cite{A2}).

The collection of eigenvalues of the  Mackey range is
\begin{align*} e(\frak T):=\{t\in\Bbb T:\ \exists\ F\in L^\infty(\frak q),\ F\nequiv 0,\ F\circ\frak T=e^{2\pi it}F\}\end{align*} and it follows from the definitions that
\begin{align*} e(\frak T)=\mathcal T(\tau,\v):=\{s\in\Bbb T:\ \exists\ f\in L^\infty(P),\ f\nequiv 0,\ ,\ f\circ\tau=e^{2\pi is\v}f\}.\end{align*}

\

It it is shown in  \S2 of \cite{AN} (see also Theorem 2.6.3 of \cite{A1}) that
\begin{align*}\tag{{\Large\Pointinghand}}\mathcal T(\tau,\v)&=G_2(b).
\end{align*}
.
\

Although formally, ({\Large\Pointinghand}) was only stated for growth sequences in \cite{AN} and \cite{A1}, the proofs do not use this condition and apply to arbitrary $b\in\incss$.

\

Consider  the Polish group  $\frak B(\mu):=\{f\in L^2(\mu):\ |f|\equiv 1\}$  equipped with $L^2(\mu)$-distance.

\

\proclaim{Lemma 4.1}
\

If the probability  $\mu\in\mathcal P(\Bbb T)$ is IP Dirichlet along $b$, then $\exists\ \mathcal X:\Om\to\frak B(\mu)$ continuous so that
\begin{align*}\tag{\dsmedical}\sup_{\om\in\Om}\,\|\chi_{b(K(\om)\cap [1,n])}-\mathcal X(\om)\|_{L^2(\mu)}\xrightarrow[n\to\infty]{}\ 0\end{align*}where $K(\om):=\{n\ge 1:\ \om_n=1\}.$\endproclaim
\demo{Proof}
\sms  Suppose that $\mu\in\mathcal P(\Bbb T)$ is IP Dirichlet along $b$. Fix $\om \in \Om$. We claim that the sequence $n\mapsto \chi_{b(K(\om)\cap [1,n])}$ is Cauchy in $L^2(\mu)$.
To see this, let
$$\mathcal E_n:=\sup_{F\in\mathcal F,\ \min\,F\ge n}\|\chi_{_{b(F)}}-1\|_{L^2(\mu)}$$
then by assumption $\mathcal E_n\xrightarrow[n\to\infty]{}\ 0$. Evidently
$$\|\chi_{b(K(\om)\cap [1,n])}-\chi_{b(K(\om)\cap [1,n+k])}\|_{L^2(\mu)}= \|\chi_{_{b(K(\om)\cap [n+1,n+k])}}-1\|_{L^2(\mu)}\le\mathcal E_n$$ whence $\exists\ \mathcal X:\Om\to \frak B(\mu)$ so that
$$\chi_{b(K(\om)\cap [1,n])}\xrightarrow[n\to\infty]{L^2(\mu)}\  \mathcal X(\om)\ \ \text{uniformly in}\ \om\in\Om.$$
For $\om\in\Om$,
$$\|\chi_{b(K(\om)\cap [1,n])}-\mathcal X(\om)\|_{L^2(\mu)}\xleftarrow[k\to\infty]{}\ \|\chi_{b(K(\om)\cap [1,n])}-\chi_{b(K(\om)\cap [1,n+k])}\|_{L^2(\mu)}\le 2\mathcal E_n$$ proving (\dsmedical). Clearly, for each $n\ge 1,\ \om\mapsto\chi_{b(K(\om)\cap [1,n])}$ is continuous ($\Om\to\frak B(\mu)$) and so continuity of $\mathcal X:\Om\to\frak B(\mu)$ follows from the
uniformity of the convergence. \ \ \Checkedbox

Note that the converse of Lemma 4.1 is also true.
\demo{Completion of the  proof }
\

Now suppose that $\mu\in\mathcal P(\Bbb T)$ is IP-Dirichlet along $b$. By Lemma 4.1,
$\exists\  \mathcal X:\Om\to\frak B(\mu)$ satisfying
(\dsmedical).

 We claim that
\begin{align*}\tag{\Radioactivity}\mathcal X(\tau\om)=\chi_{\varphi(\om)}\mathcal X(\om).\end{align*}
To see this, note that $\chi_{b(K(\om)\cap [1,n])}(t)=\prod_{k=1}^n\chi_{\om_kb_k}(t)=:\mathcal X_n(\om,t)$. For $n>\ell(\om)$,
$$\frac{\mathcal X_n(\tau\om)}{\mathcal X_n(\om)}=\frac{\chi_{b(K(\tau\om)\cap [1,n])}}{\chi_{b(K(\om)\cap [1,n])}}=
\prod_{k=1}^n\frac{\chi_{(\tau\om)_kb_k}}{\chi_{\om_kb_k}}=\chi_{\varphi(\om)}.$$
Since $\frac{\mathcal X_n(\tau\om)}{\mathcal X_n(\om)}\xyr[n\to\infty]{\frak B(\mu)}\frac{\mathcal X(\tau\om)}{\mathcal X(\om)}$, this proves (\Radioactivity).
\

By (\dsmedical), $\exists\ n_J\to\infty$ so that
$$\sum_{J\ge 1}\|\mathcal X_{n_J}-\mathcal X\|_{L^2(P\x\mu)}<\infty\le \sum_{J\ge 1}\sup_{\om\in\Om}\,\|\chi_{b(K(\om)\cap [1,n_J])}-\mathcal X(\om)\|_{L^2(\mu)}<\infty$$ (where $\mathcal X(\om,t):=\mathcal X(\om)(t)$). Hence $\mathcal X_{n_J}\to\mathcal X$
$P\x\mu$-a.e. and by Fubini's theorem,
 $\exists\ \Lambda\in\B(\Bbb T),\ \mu(\Lambda)=1$ so that
$$\prod_{k=1}^{n_J}\chi_{\om_kb_k}(t)\xrightarrow[J\to\infty]{}\ \mathfrak X_t(\om)=\mathcal X(\om)(t)\ \forall\ t\in\Lambda\ \&\ P-\text{a.e.}\ \om\in\Om.$$
By (\Radioactivity), for $t\in\Lambda$,
$$\mathfrak X_t\circ\tau=e^{2\pi it\v}\mathfrak X_t\ \ P-\text{a.e.}$$
and $t\in \mathcal T(\tau,\v)$. By ({\Large\Pointinghand}), $t\in G_2(b)$.\ \ \ \ \Checkedbox

The following example shows that the converse to Theorem 2 is false.
\proclaim{Example 4.2}\ \ \ \ $\exists\ b\in\Bbb N^\Bbb N(\uparrow)\ \&\ \mu\in\mathcal P(\Bbb T)$ non-atomic, not IP-Dirichlet along $b$ but so that $\mu(G_2(b))=1$.\endproclaim
\demo{\tt Construction:}
\

Define $b\in\Bbb N^\Bbb N(\uparrow)$ by $b_n:=\prod_{k=1}^na_k$ with $a_k:=k+1$.
\

Consider the mapping
$t:\Om:=\prod_{k\ge 1}\{0,1\}\to [0,1]$ defined by
$$t(\om)=t(\om_1,\om_2,\dots):=\sum_{n\ge 1}\frac{\om_n}{b_n}.$$
This is injective and Borel measurable, so $t(\Om)$ is an uncountable, Borel set in $[0,1]$. We claim that
\begin{align*}\tag{\Biohazard}t(\Om)\subset G_2(b).\end{align*}
\demo{Proof of (\Biohazard)}

For $\om\in\Om$  and $N\ge 1$, we have that
$$b_Nt(\om)=\frac{\om_{N+1}}{a_{N+1}}+\frac{1}{a_{N+1}}\sum_{k\ge 2}\frac{\om_{N+k}}{a_{N+2}\dots a_{N+k}}\ \ \ \mod 1.$$

Now
$$|\sum_{k\ge 2}\frac{\om_{N+k}}{a_{N+2}\dots a_{N+k}}|\le\sum_{k\ge 2}\frac1{a_{N+2}\dots a_{N+k}}<\frac1{a_{N+1}}.$$
Thus,  we have
\begin{align*}\tag{\Laserbeam}\<b_Nt(\om)\>=\frac{\om_{N+1}}{a_{N+1}}+\th_N(\om)\ \ \text{where}\ \ \th_N\le \tfrac{1}{a_{N+1}^2},\end{align*}
   whence $|\<b_Nt(\om)\>|^2\le \tfrac2{N^2}$ and
$$\sum_{k=1}^\infty|\<b_kt(\om)\>|^2\le\frac{\pi}3.\ \ \ \text{\Checkedbox}\ \ \text{(\Biohazard)}$$

\

Now define $P\in\mathcal P(\Om)$ by $P:=\prod_{k\ge 1}(\tfrac12\d_0+\tfrac12\d_1)$ and set $\mu:=P\circ{t}^{-1}$, then $\mu\in\mathcal P(G_2(b))$. We now show that $\mu$ is not IP-Dirichlet along $b$.

\

Fix $1<\l<e^{\frac13}$, then by (\Laserbeam),
$$\big|\sum_{\l^N<k<\l^{N+1}}\<b_Nt(\om)\>-\sum_{\l^N<k<\l^{N+1}}\frac{\om_{N+1}}{a_{N+1}}\big|\le \sum_{\l^N<k<\l^{N+1}}\tfrac2{k^2}\xrightarrow[N\to\infty]{}0$$
and setting $s_n:=\sum_{j=1}^n\om_j, \kappa_n:=\lcl\l^n\rcl,\ \ell_n:=\lfl\l^{n+1}\rfl$, we have
\begin{align*}\sum_{\l^N<k<\l^{N+1}}\frac{\om_{k}}{a_{k}}&=\sum_{\kappa_N\le k\le \ell_N}\frac{s_k-s_{k-1}}{k+1}\\ &=\sum_{\kappa_N\le k\le \ell_N}\frac{s_k}{k+1}-\sum_{\kappa_N-1\le k\le \ell_N-1}\frac{s_k}{k+2}\\ &=
\frac{s_{\ell_N}}{\ell_N+1}-\frac{s_{\kappa_N}}{\kappa_N+2}+\sum_{\kappa_N\le k\le \ell_N-1}\frac{s_k}{(k+1)(k+2)}.\end{align*}
By the {\tt SLLN}, we have
$$\frac{s_N}{N}\xrightarrow[N\to\infty]{}\frac12\ \ \ \text{a.s.}$$
Writing $\G_N\approx\D_N$ as $N\to\infty$ to mean $\G_N\approx\D_N\xyr[N\to\infty]{}\ 0$, we have
 for $P$-a.e. $\om\in\Om$, as $N\to\infty$:
\begin{align*}\sum_{\l^N<k<\l^{N+1}}\<b_Nt(\om)\> &\approx \sum_{\l^N<k<\l^{N+1}}\frac{\om_{k}}{a_{k}}\\ &=
\frac{s_{\ell_N}}{\ell_N+1}-\frac{s_{\kappa_N}}{\kappa_N+2}+\sum_{\kappa_N\le k\le \ell_N-1}\frac{s_k}{(k+1)(k+2)}\\ &\approx\sum_{\kappa_N\le k\le \ell_N-1}\frac1{2k}\\ &\lra \log\l.\end{align*}
Thus
\begin{align*}\sup_{F\in\mathcal F,\ \min\,F>\l^N}\|\chi_{b(F)}-1\|_{L^2(\mu)}^2 &\ge
\Bbb E(|\exp[2\pi i\sum_{\l^N<k<\l^{N+1}}b_Nt(\om)]-1|^2)\\ &\ge
4\Bbb E(|\sum_{\l^N<k<\l^{N+1}}\<b_Nt(\om)\>|^2)\\ &\lra 2\log\l.\ \ \ \boxtimes \end{align*}

\

\section*{\S5  Remarks on the thickness of rigidity sequences}
\

\subsection*{Remark 5.1} \ \  Rigidity sequences for weakly mixing transformations can be arbitrarily "large" within the limitation of density zero.
\

It follows from the definitions that for $b\in\incss$ a  growth sequence,
$$|{\tt FS}\,(b)\cap [1,n]|\asymp 2^{c(n)}$$ where $c(n)=\min\,\{k\ge 1:\ b_k\ge n\}$ and if $T\in {\tt MPT}\xbm$ is IP-rigid along $b$, one might expect a rigid sequence at least of this thickness.
\

Indeed, in this case,  by Theorem 2, $\s_T(e(\frak T)^c)=0$ where $\frak T$ is the dyadic tower with growth sequence $b$ (since by [AN] $e(\frak T)=G_2(b)$). By theorem 4 in \cite{A2} $\exists\ L\subset\Bbb N$ with
$$T^n\xyr[n\to\infty,\ n\in L]{{\tt MPT}}\ \text{\tt Id}\ \ \&\ \ \frac{|L\cap [1,n]|}{2^{c(n)}}\xyr[n\to\infty]{}\infty.$$
\

By the Corollary in \cite{A2}, $\forall\ a(n)>0,\ \tfrac{a(n)}n\xyr[n\to\infty]{}0$, there is a    weakly mixing $T\in{\tt MPT}\xbm$ and
 $L\subset\Bbb N$ such that
 $$T^n\xyr[n\to\infty,\ n\in L]{{\tt MPT}}\ \text{\tt Id}\ \ \&\ \ \frac{|L\cap [1,n]|}{a(n)}\xyr[n\to\infty]{}\infty.$$
For more on this phenomenon, see \S3 in \cite{BJLR}.
\

\subsection*{Remark 5.2}\ \  Let $b\in\incss$ be a growth sequence and suppose that ${\tt H-dim}(G_1(b))>\a$ (where $\a\in (0,1)$ and {\tt H-dim} denotes Hausdorff dimension).
\

We claim that  $\exists\ T\in{\tt MPT}$ weakly mixing $\&$ IP-rigid along $b$ with the property that for any sequence $L\subset\Bbb N$ along which $T$ is rigid:
\begin{align*}\tag{\Football}\sum_{n=1}^\infty\frac{|L\cap [1,n]|}{n^{2-\a}}<\infty.\end{align*} Note that it follows from this that
$\sum_{n=1}^\infty\frac{2^{c(n)}}{n^{2-\a}}<\infty.$\demo{Proof of (\Football)}
\

As in the proof of theorem 1 of \cite{A3}, it follows from Frostman's theorem (\cite{Fr}, see also \cite{KS}) that  $\exists\ \mu\in\mathcal P(G_1(b))$  so that
\begin{align*} \sum_{n=1}^\infty\frac{|\widehat{\mu}(n)|}{n^{1-\a}}<\infty.\end{align*}
Let $T\in{\tt MPT}\xbm$ be the associated Gaussian automorphism. By Proposition 1, $T$ is IP-rigid along $b$.
\

Suppose that $T$ is rigid along $L\subset\Bbb N$. In particular $\exists\ N\ge 1$ so that $|\widehat{\mu}(n)|\ge\tfrac12\ \forall\ n\in L,\ n>N$.
\

It follows that
\begin{align*}\sum_{n=1}^\infty\frac{|L\cap [1,n]|}{n^{2-\a}}&=
\sum_{n=1}^\infty \sum_{k=1}^n1_L(k)\frac1{n^{2-\a}}\\ &=
\sum_{k=1}^\infty 1_L(k)\sum_{n=k}^\infty\frac1{n^{2-\a}}\\ &\le
\frac1{1-\a}\sum_{k=1}^\infty 1_L(k)\frac1{k^{1-\a}}\\ &\le\frac{N}{1-\a}+
\frac2{1-\a}\sum_{k=N+1}^\infty \frac{|\widehat{\mu}(k)|}{k^{1-\a}}\\ &<\infty.\ \ \ \ \text{\Checkedbox}\text{(\Football)}\end{align*}
\subsection*{Remark 5.3}\ \   The condition (\Football) is sharp. In \S5 of \cite{A3}, $\forall\ \a\in (0,1)$ a growth sequence $b^{(\a)}\in\incss$ is exhibited with
$$\text{\tt H-dim}\,(G_1(b^{(\a)}))=\a\ \&\ 2^{c^{(\a)}(n)}\gg n^{1-\a}$$
whence
 if $T\in {\tt MPT}\xbm$ is so that $\s_T(G_2(b))=1$, then (again by theorem 4 in [A2]) $\exists\ L\subset\Bbb N$ rigid for $T$ with $\frac{|L\cap [1,n]|}{n^{1-\a}}\to\infty$ and therefore
 \begin{align*}\sum_{n=1}^\infty\frac{|L\cap [1,n]|}{n^{2-\a}}=\infty.\end{align*}

\end{document}